\documentclass[a4paper]{article}

\usepackage{amsmath,amssymb}
\usepackage[dvips]{epsfig}
\usepackage[latin1]{inputenc}

\newcommand{\ZZ}{\mathbb{Z}}
\newcommand{\xn}{{X_n}}
\newcommand{\hh}{\tilde{h}}
\newcommand{\uu}{\tilde{u}}
\newcommand{\gr}{\operatorname{gr}}
\newcommand{\ach}{\textsf{A}_\xn}
\newcommand{\fil}{\mathcal{F}}
\newcommand{\phip}{\Phi_{+}}

\newtheorem{theo}{Theorem}[section]
\newtheorem{prop}[theo]{Proposition}
\newtheorem{lemma}[theo]{Lemma}

\newenvironment{proof}{\begin{trivlist}\item{\bf{Proof.}}}
  {\hfill\rule{2mm}{2mm}\end{trivlist}}

\title{Antichains of positive roots and Heaviside functions}
\date{\today} \author{Frédéric Chapoton}

\begin{document}
\maketitle
\begin{abstract}
  The ring of locally-constant integer-valued functions on the
  dominant chamber of the Shi arrangement is endowed with a filtration
  and a new basis, compatible with this filtration, is found. This
  basis is compared to the trivial basis. The ring is given a
  presentation by generators and relations.
\end{abstract}

\setcounter{section}{-1}
\section{Introduction}

The aim of this short paper is to study some rings associated to
finite root systems. These rings are defined starting from a
well-known hyperplane arrangement associated to a root system, which
is called the Shi arrangement \cite{shi1,shi2}. The number of connected
components of the Shi arrangement which are contained in the dominant
Weyl chamber is known to be the number of antichains of the poset of
positive roots for the standard order and is equal to the generalized
Catalan number which also appears in the theory of cluster algebras of
Fomin and Zelevinsky \cite{fomzel1,fomzel2}, see
\cite{athanas1,athanas2,cellini} and references therein. The main
object in the present article is the ring of locally-constant
integer-valued functions on the intersection of the dominant Weyl
chamber with the Shi arrangement. This commutative ring, which is a
free abelian group of rank given by the generalized Catalan number
above, is endowed with a filtration by a general construction on
hyperplane arrangements due to Gelfand and Varchenko \cite{gelfvarch}.
A presentation by generators and relations is given, which leads to a
new basis indexed by antichains and compatible with the filtration.
The relation between the basis of idempotents and the new basis is
explained using a natural partial order on antichains.

\section{The root order}

Let $\xn$ be a finite Dynkin diagram, \textit{i.e.} $X_n$ is either in
the classical series $A_n,B_n,C_n,D_n$ or one of the exceptionals
$E_6,E_7,E_8,F_4,G_2$.

Let $\phip$ be the set of positive roots for $\xn$. This set is
endowed with the standard partial order $\leq$ defined by $\alpha \leq
\beta$ if the difference $\beta-\alpha$ has non-negative coefficients
in the basis of simple roots.

Let $\ach$ be the set of antichains in the poset $(\phip,\leq)$. The
elements of $\ach$ are called non-nesting partitions in the literature
\cite{athanas1,cellini,reiner}.

In any poset, there is a simple bijection between antichains and upper
ideals. An antichain is mapped to the upper ideal it generates and an
upper ideal is mapped to the set of its minimal elements, which is an
antichain. The upper ideal associated to an antichain $p$ is denoted
by $I_p$.

\section{The Shi hyperplane arrangements}

The roots in $\phip$ are considered as linear forms on a real vector
space of dimension $n$ in the usual way.

The Shi hyperplane arrangement associated to $\xn$ is the collection
formed by all the hyperplanes $\alpha(x)=0$ and $\alpha(x)=1$ where
$\alpha$ describes $\phip$.

Let us call \textit{region} a connected component of the complement of
the union of these hyperplanes and \textit{dominant region} a region
contained in the dominant Weyl chamber.

By works of Athanasiadis and Cellini-Papi (see \cite{athanas2} and
references therein), the number of dominant regions is known to be the
generalized Catalan number associated to $X_n$, which is
\begin{equation}
  C_\xn=\prod_{i=1}^{n}\frac{h+e_i+1}{e_i+1},
\end{equation}
where $h$ is the Coxeter number and $e_1,\dots,e_n$ the exponents of
$\xn$.

Let us recall here the bijection between dominant regions and
antichains. An antichain $p$ is mapped to the dominant region defined
by $\alpha(x)>1$ for all $\alpha \in I_p$ and $0<\alpha(x)<1$
elsewhere. One recovers the ideal $I_p$ from a dominant region as the
set of roots which take values greater than $1$ on this region.

Define a partial order $\preceq$ on antichains by inclusion of
associated upper ideals, \textit{i.e.} set $p \preceq q$ if $I_p
\subseteq I_q$ as a set. This is just the lattice of upper ideals of
$(\phip,\leq)$ for the inclusion order.

\begin{figure}
  \begin{center}
    \leavevmode 
    \epsfig{file=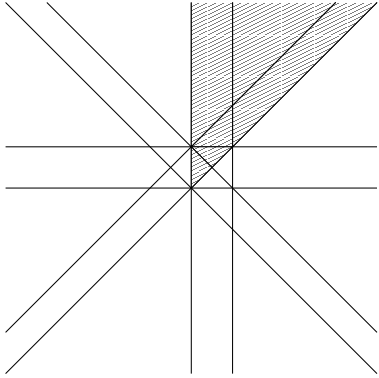,width=5cm} 
    \caption{Shi arrangement for $B_2$}
    \label{chambre}
  \end{center}
\end{figure}

Fig. \ref{chambre} displays the Shi arrangement and the dominant Weyl
chamber for $B_2$.

\section{Heaviside functions and filtration}

Let $V_{\xn}$ be the ring of locally-constant functions with integer
values on the complement of the Shi arrangement of $\xn$. 

For any hyperplane arrangement, an increasing filtration on the ring
of locally-constant integer-valued functions has been introduced by Gelfand
and Varchenko in \cite{gelfvarch}. In the case of the Shi arrangement,
this filtration is defined as follows. First for each root $\alpha \in
\phip$, define two locally-constant functions $h^0_\alpha$ and
$h^1_\alpha$ as follows.
\begin{equation}
  h^0_\alpha(x)=
  \begin{cases}
    0 \text{ if }\alpha(x)<0, \\
    1 \text{ if }\alpha(x)>0,
  \end{cases}
\text{ and  } 
  h^1_\alpha(x)=
  \begin{cases}
    0 \text{ if }\alpha(x)<1, \\
    1 \text{ if }\alpha(x)>1.
  \end{cases}
\end{equation}
These are called Heaviside functions by Gelfand and Varchenko
\cite{gelfvarch} by similarity with the step-function of Heaviside.

Then the filtration $\fil$ is
\begin{equation}
  \ZZ 1 = \fil_0\subseteq \fil_1
  \subseteq \fil_2 \subseteq \dots \subseteq \fil_n \subseteq \dots ,
\end{equation}
where $\fil_1$ is the linear span of the functions $1$, $h^0_\alpha$,
$h^1_\alpha$ for all positive roots $\alpha$ and the space $\fil_k$ is
defined to be $(\fil_1)^k$ for all $k>1$.

A key result of Gelfand and Varchenko in \cite{gelfvarch} says that
the filtration reaches the full ring at step $n$, \textit{i.e.}
$\fil_n=V_\xn$. This implies that $V_\xn$ is generated by the
Heaviside functions $h^0_\alpha$, $h^1_\alpha$.

\section{Restriction to the dominant chamber}

Let $H_{\xn}$ be the ring of locally-constant functions with integer
values on the intersection of the complement of the Shi arrangement
and the dominant Weyl chamber. The rank of $H_\xn$ is therefore the
number $C_\xn$ of dominant regions of the Shi arrangement.

Let $p$ be an antichain. By the correspondence between antichains and
dominant regions, one can define a function $\delta_p$ in $H_\xn$
which has value $1$ on the region corresponding to $p$ and vanishes
elsewhere.

The set of functions $(\delta_p)_p$ where $p$ describes $\ach$ is of
course a basis of $H_\xn$ made of orthogonal idempotents, called the
trivial basis.

As the set of dominant regions is a subset of the set of regions,
there is a surjective restriction morphism from $V_{\xn}$ to $H_\xn$.
The filtration $\fil$ of $V_\xn$ induces a filtration on $H_\xn$ still
denoted by $\fil$. Remark that the image of $h^0_{\alpha}$ is $1$ for
all $\alpha$. Denote simply by $h_\alpha$ the image of $h^1_\alpha$.
So the filtration $\fil$ on $H_\xn$ is given by
\begin{equation}
  \ZZ 1 = \fil_0\subseteq \fil_1
  \subseteq \fil_2 \subseteq \dots \subseteq \fil_n \subseteq \dots ,
\end{equation}
where $\fil_1$ is the linear span of the functions $1$,$h_\alpha$ for
all positive roots $\alpha$ and the space $\fil_k$ is
defined to be $(\fil_1)^k$ for all $k>1$.

From surjectivity and the similar result for $V_\xn$, one has
$\fil_n=H_\xn$. In particular $H_\xn$ is generated by the functions
$h_\alpha$.

\begin{lemma}
  \label{hachedel}
 In the basis $(\delta_p)_p$ of $H_\xn$, the function $h_\alpha$ is
  given by
  \begin{equation}
    h_\alpha=
    \sum_{\{\alpha\} \preceq p} \delta_{p}.
  \end{equation}
\end{lemma}
\begin{proof}
  From its definition by restriction of a Heaviside function,
  $h_\alpha$ is the characteristic function of the regions where
  $\alpha$ takes values greater than $1$. Therefore it is the sum of
  $\delta_p$ over each dominant region $p$ where $\alpha$ takes values
  greater than $1$. This condition means exactly that $\alpha \in I_p$
  or $I_{\{ \alpha \}} \subseteq I_p$, \textit{i.e.} $\{\alpha\}
  \preceq p$.
\end{proof}

\begin{lemma}
  The map $\alpha \mapsto \{ \alpha\}$ gives an order reversing
  injection of the poset of roots $(\phip,\leq)$ in the poset of
  antichains $(\ach,\preceq)$.
\end{lemma}
\begin{proof}
  Obvious.
\end{proof}

\begin{prop}
  The functions $h_\alpha$ satisfy
\begin{equation}
  \label{relah}
  h_\alpha h_\beta = h_{\min(\alpha,\beta)},
\end{equation}
whenever $\{\alpha,\beta\}$ is not an antichain.
\end{prop}
\begin{proof}
  Assume for example that $\alpha \leq \beta$, so that
  $\min(\alpha,\beta)=\alpha$. Then $I_{\{\beta\}}\subseteq
  I_{\{\alpha\}}$ and $\beta \preceq \alpha$. Therefore using the
  description of $h_\alpha$ as a sum of idempotents given in Lemma
  \ref{hachedel}, one has $h_{\alpha}h_{\beta}=h_{\alpha}$.
\end{proof}

%\section{Associated graded algebra}

Let $\gr H_\xn$ be the graded ring associated to the filtration $\fil$
of the ring $H_\xn$. Then $\gr H_\xn$ is generated by elements
$\hh_\alpha$ which satisfy
\begin{equation}
  \hh_\alpha \hh_\beta =0,
\end{equation}
whenever $\{\alpha,\beta\}$ is not an antichain.

\section{Rings presented by generators and relations}

Consider the commutative ring $U_\xn$ generated by the variables
$u_\alpha$ for $\alpha \in \phip$ modulo the relations
\begin{equation}
  \label{relau}
  u_\alpha u_\beta = u_{\min(\alpha,\beta)},
\end{equation}
whenever $\{\alpha,\beta\}$ is not an antichain.

One can define an increasing filtration $\fil$ on $U$ as follows
\begin{equation}
    \ZZ 1 = \fil_0\subseteq \fil_1
  \subseteq \fil_2 \subseteq \dots \subseteq \fil_n \subseteq \dots ,
\end{equation}
where $\fil_1$ is the linear span of the elements $1$,$u_\alpha$ for
all positive roots $\alpha$ and the space $\fil_k$ is defined to be
$(\fil_1)^k$ for all $k>1$.

Let $\gr U_\xn$ be the associated graded ring. Then $\gr U_\xn$ is
presented by the generators $\uu_\alpha$ modulo the relations
\begin{equation}
  \label{relauu}
  \uu_\alpha \uu_\beta = 0,
\end{equation}
whenever $\{\alpha,\beta\}$ is not an antichain.

Let us define an element $u_p$ of $U_\xn$ for each antichain $p$ as
follows. If $p=\{\alpha_1,\dots,\alpha_k\}$ then set
$u_p=u_{\alpha_{1}} \dots u_{\alpha_{k}}$. Note that there is no
conflict in notation. In the same way, one defines elements $\uu_{p}$
in $\gr U_\xn$.

\begin{prop}
  \label{monomebase}
  The set $(\uu_p)_p$ is a basis of the ring $\gr U_\xn$. 
\end{prop}
\begin{proof}
  This set of monomials spans $\gr U_\xn$ because any monomial
  containing two comparable roots vanish by relations (\ref{relauu}).
  Now it follows immediately from the shape of relations (\ref{relauu})
  that the ring $\gr U_\xn$ is in fact graded by the free abelian
  group generated by the set $\phip$. So each of the monomials $\uu_p$
  has a different weight. Hence there can be no linear relation
  between these monomials except maybe that some of them may vanish.
  But one can see that monomials which vanish because of relations
  (\ref{relauu}) are precisely the monomials containing two comparable
  roots.
\end{proof}

\begin{prop}
  The set $(u_p)_p$ is a basis of the ring $U_\xn$. 
\end{prop}
\begin{proof}
  This is a direct corollary of Theorem \ref{monomebase}.
\end{proof}

Therefore the rank of $U_\xn$ is the generalized Catalan number
$C_\xn$. Furthermore the rank of the graded component of degree $k$ of
$\gr U_\xn$ is the number of antichains of cardinal $k$.

\section{Isomorphism}

\begin{prop}
  There exists a unique morphism of filtered rings $\rho$ from $U_\xn$
  to $H_\xn$ which sends $u_\alpha$ to $h_\alpha$.
\end{prop}
\begin{proof}
  It is enough to check that the relations (\ref{relau}) are satisfied
  in $H_\xn$. But this is exactly (\ref{relah}). The filtrations are
  clearly mapped one to another.
\end{proof}

\begin{prop}
  The morphism $\rho$ is surjective.
\end{prop}
\begin{proof}
  The image of $\rho$ contains the elements $h_\alpha$ which generate $H_\xn$.
\end{proof}

\begin{theo}
  \label{principal}
  The morphism $\rho$ is an isomorphism of filtered rings between
  $U_\xn$ and $H_\xn$.
\end{theo}
\begin{proof}
  The morphism $\rho$ preserves filtrations, is surjective and both
  rings have the same rank given by the generalized Catalan number
  $C_\xn$.
\end{proof}

\begin{prop}
  The induced morphism $\rho$ from $\gr U_\xn$
  to $gr H_\xn$ is an isomorphism of graded rings.
\end{prop}
\begin{proof}
  This is a corollary of Theorem \ref{principal}.
\end{proof}

\section{Change of basis}

The aim of this section is to study the relation between the bases
$(h_p)_p$ and $(\delta_p)_p$ of $H_\xn$, both indexed by the set
$\ach$ of antichains.

\begin{prop}
  Let $p$ be an antichain. One has
  \begin{equation}
      h_p= \sum_{p \preceq q} \delta_q.
  \end{equation}
\end{prop}
\begin{proof}
  By induction on the cardinal of the antichain. This works for the
  empty antichain. This is also true if the antichain is a singleton
  by Lemma \ref{hachedel}. Assume that it is proven for antichains
  with less elements. Let $p=p'\sqcup \{\alpha\}$, so that
  $h_{p'}=h_{p}h_{\alpha}$. Then
  \begin{equation}
    h_{p}=\sum_{p'\preceq q} \sum_{\{\alpha \} \preceq r}\delta_{q} \delta_{r}.
  \end{equation}
  From the idempotency and orthogonality of the basis $(\delta_p)_p$,
  one has
  \begin{equation}
    h_{p}=\sum_{(p'\sqcup \{\alpha \}) \preceq q} \delta_{q},
  \end{equation}
  because the union (which is also the supremum) of the ideals
  $I_{p'}$ and $I_{\{\alpha\}}$ is the ideal $I_{p'\sqcup \{\alpha
    \}}=I_p$.
\end{proof}

So the coefficient matrix of the basis $(h_p)_p$ in the basis
$(\delta_p)_p$ is given by the poset matrix for $\preceq$.

Hence, conversely, the coefficients of the idempotents $\delta_p$ in
the basis $(h_p)_p$ are described by the Möbius matrix of the poset
$(\ach,\preceq)$.

\begin{figure}
  \begin{center}
    \leavevmode 
    \epsfig{file=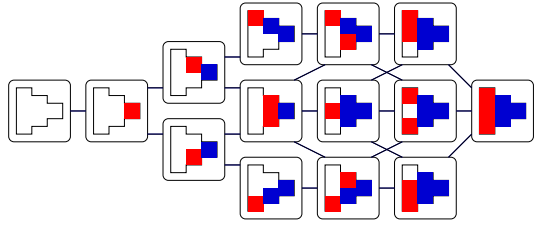,width=7cm} 
    \caption{Poset of antichains as upper ideals}
    \label{inclusion}
  \end{center}
\end{figure}

Fig. (\ref{inclusion}) displays the poset of antichains for the Dynkin
diagram $A_3$.

\smallskip

Let us remark that there is a simple bijection between antichains of
type $A_n$ and Dyck paths which maps the cardinal of the upper ideal
associated to an antichain to the area above the Dyck path. Therefore
the antichains of $\xn$ counted according to the cardinal of the upper
ideal give a possible generalization to root systems of the classical
$q$-Catalan numbers corresponding to Dyck paths and area
\cite{carlitz}.

\nocite{*}

\bibliographystyle{plain}
\bibliography{antich}

\begin{thebibliography}{10}

\bibitem{athanas1}
Christos~A. Athanasiadis.
\newblock On noncrossing and nonnesting partitions for classical reflection
  groups.
\newblock {\em Electron. J. Combin.}, 5(1):Research Paper 42, 16 pp.
  (electronic), 1998.

\bibitem{athanas2}
Christos~A. Athanasiadis.
\newblock Generalized {C}atalan numbers, {W}eyl groups and arrangements of
  hyperplanes.
\newblock preprint, 2002.

\bibitem{carlitz}
L.~Carlitz and J.~Riordan.
\newblock Two element lattice permutation numbers and their
  {$q$}-generalization.
\newblock {\em Duke Math. J.}, 31:371--388, 1964.

\bibitem{cellini}
P.~Cellini and P.~Papi.
\newblock ad-nilpotent ideals of a {B}orel subalgebra {II}.
\newblock {\em J. Algebra}, (258):112--121, 2002.

\bibitem{fomzel1}
Sergey Fomin and Andrei Zelevinsky.
\newblock Cluster algebras. {I}. {F}oundations.
\newblock {\em J. Amer. Math. Soc.}, 15(2):497--529 (electronic), 2002.

\bibitem{fomzel2}
Sergey Fomin and Andrei Zelevinsky.
\newblock {Cluster algebras II: Finite type classification}.
\newblock {\em Inventiones Mathematicae}, 2003.
\newblock {arXiv:math.RA/0208229}.

\bibitem{reiner}
Victor Reiner.
\newblock Non-crossing partitions for classical reflection groups.
\newblock {\em Discrete Math.}, 177(1-3):195--222, 1997.

\bibitem{shi1}
Jian~Yi Shi.
\newblock Sign types corresponding to an affine {W}eyl group.
\newblock {\em J. London Math. Soc. (2)}, 35(1):56--74, 1987.

\bibitem{shi2}
Jian-Yi Shi.
\newblock The number of {$\oplus$}-sign types.
\newblock {\em Quart. J. Math. Oxford Ser. (2)}, 48(189):93--105, 1997.

\bibitem{gelfvarch}
A.~N. Varchenko and I.~M. Gelfand.
\newblock Heaviside functions of a configuration of hyperplanes.
\newblock {\em Funktsional. Anal. i Prilozhen.}, 21(4):1--18, 96, 1987.

\end{thebibliography}

\end{document}